\documentclass[10pt]{article}
\usepackage{amssymb, amsmath}
\usepackage{graphicx}

\begin{document}
\title{How many times can the volume of a convex\\ 
polyhedron be increased by isometric deformations?}
\author{Victor Alexandrov}
\date{March 1, 2017}
\maketitle
\begin{abstract}
We prove that the answer to the question of the title is
`as many times as you want.' More precisely, given any constant $c>0$,
we construct two oblique triangular bipyramids, $P$ and $Q$, such that 
$P$ is convex,  
$Q$ is nonconvex and intrinsically isometric to $P$, and  
$\textrm{vol\,}Q>c\cdot \textrm{vol\,}P>0$.
\end{abstract}

\textbf{1. Introduction.}
According to the classical theorem by A.L. Cauchy and A.D. Alexandrov, 
two compact boundary-free convex polyhedral surfaces in Euclidean 3-space
are necessarily congruent as soon as they are isometric in their intrinsic metrics, 
see, e.\,g., \cite{Al05}.
Obviously, this is not true if at least one of the surfaces is nonconvex.
In \cite{BZ95}, the authors
studied isometric immersions of polyhedral surfaces and, among other things,
proved that there exists a compact 
boundary-free convex polyhedral surface allowing another isometric immersion which,
being a nonconvex polyhedral surface, encloses a larger volume than that enclosed
by the original convex surface.
This amazing existence theorem has gained new significance after 
the remarkable contribution of D.D. Bleecker, who explicitly built volume increasing isometric 
deformations of the surfaces of the Platonic solids, 
see \cite{Bl96}.
For example, he has shown that the surface 
of a regular tetrahedron can be isometrically 
deformed in such a way as to enlarge the enclosed volume by 37.7\%. 
The results of D.D. Bleecker were presented in popular literature, see 
\cite{Al00}  and \cite{Pa08}.
They were also improved in \cite{MG09}  
by constructing more sophisticated isometric deformations of the surfaces
of the Platonic solids; for example, it was shown that the surface of a regular tetrahedron 
can be isometrically deformed in such a way as to enlarge the enclosed volume by 44\%.
The next contribution to this field of geometry was made in
\cite{Sa10},
where it is proved that every compact boundary-free convex polyhedral surface 
in Euclidean 3-space possesses a volume increasing isometric deformation.

The main result of this paper is the following 

\textbf{Theorem 1.}\textit{Given any constant $c>0$,
there exist two oblique triangular bipyramids, $P$ and $Q$, 
in Euclidean 3-space such that}

(i) $P$ \textit{is convex and encloses a non-zero volume
{\rm (}i.\,e., $\text{{\rm vol\,}}P>0${\rm )},}  

(ii)  $Q$ textit{is nonconvex and intrinsically isometric to $P$, and}

(iii)  $\text{{\rm vol\,}}Q>c\cdot \text{\rm vol\,}P$.

\textit{Remark 1.} 
In all the articles mentioned above, the deformed nonconvex surfaces have much more 
faces then the original convex surface, e.\,g, at least 7 times more in \cite{Bl96}.
In contrast, in Theorem~1 the numbers of vertices, edges, and faces are the same for the convex
bipyramid $P$ and its nonconvex deformed counterpart~$Q$.

\textit{Remark 2.} 
It follows immediately from the classical theorem by A.L. Cauchy and A.D. Alexandrov
that the bipyramids $P$ and $Q$ from Theorem~1 cannot be included in a continuous
family of pairwise intrinsically isometric polyhedral surfaces.
Hence, Theorem~1 does not belong to the theory of flexible polyhedra. 
For more details about this theory the reader is referred to 
\cite{Sa11}, \cite{Ga14}, and the references given there.

\textbf{2. Basic definitions.}
In this section, we specify the terminology used.

We study compact boundary-free polyhedral surfaces in Euclidean 3-space.
For short, we call them polyhedral surfaces or polyhedra.
A polyhedral surface is called convex if it is the boundary of a convex set.
Bipyramid is any polyhedral surface combinatorially equivalent
to the boundary of the convex hull of a regular tetrahedron and its
image under reflection in one of its faces.
Two polyhedral surfaces, $p$ and $q$, are called intrinsically isometric 
if there exists a one-to-one correspondence $f:p\to q$,
which preserves the length of any curve.
Such $f:p\to q$ is called an intrinsic deformation of a polyhedral surface $p$.
An intrinsic deformation $f:p\to q$ is called volume increasing if $p$ is 
orientable and the 3-dimensional volume enclosed by $p$
is less than the 3-dimensional volume enclosed by $q$.

\textbf{3. Auxiliary polyhedron $p(t)$.}
In this section, we build and study a continuous family of auxiliary 
convex oblique triangular bipyramids $p(t)$. 
The parameter $t$ is a real number, which takes values from
an interval $0<t<t_0$.
The exact value of $t_0$ will be specified below.
\begin{figure}
\begin{center}
\includegraphics[width=0.4\textwidth]{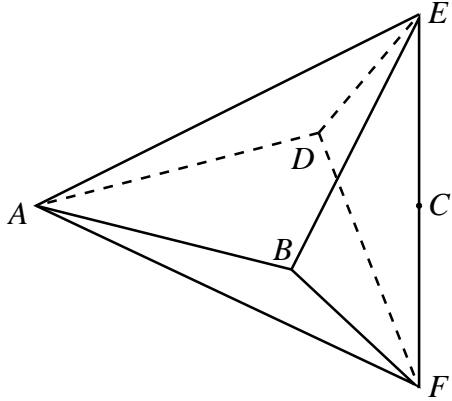}
\end{center}
\caption{Convex oblique triangular bipyramid $p(t)$}
\end{figure}

The convex oblique triangular bipyramid $p(t)$ is shown schematically in Fig.~1.
The point $C$ is not a vertex of $p(t)$. $C$ is the middle point of the edge $EF$.
It is shown in Fig.~1 for convenience of further explanations.
The bipyramid $p(t)$ is assumed to be symmetric with respect to the planes of the triangles
$ABD$ and $AEF$.
By definition, the lengths of some edges of $p(t)$ are as follows:
$|AB|=|AD|=10$, $|BE|=|BF|=|DE|=|DF|=13$, and $|EF|=24$.
The lengths of the edges $AE$ and $AF$ are supposed to be equal to each other
and depend on the parameter $t$, which has the following geometric meaning:
$t$ is equal to the nonoriented angle between the straight line segments 
$AB$ and $AC$. According to this definition, $t$ may take any value
from 0 to $t_0=\pi/6$. The value $t_0=\pi/6$ corresponds to the case, when
the triangles $BEF$ and $DEF$ lie in a common plane and, thus, the
quadrilateral $ABCD$ is in fact the regular triangle, whose sides
have the lengths 10.

In order to express $|AE|$ as a function of $t$, we proceed as follows.
Applying the Pythagorean theorem to the right triangle $BCE$, we find $|BC|=5$.
The cosine law applied to the triangle $ABC$ yields
$|BC|^2=|AB|^2+|AC|^2-2|AB||AC|\cos t$. Thus,
$|AC|=10\cos t\pm\sqrt{100\cos^2t-75}.$
In order to obtain a convex bipyramid $p(t)$, 
we must choose the maximal value of $|AC|$ in the previous formula. 
So, we get
$|AC|=10\cos t+\sqrt{100\cos^2t-75}$.
Now, applying the Pythagorean theorem to the right triangle $ACE$, we find
\begin{equation}\label{eqn:1}
|AE|=\sqrt{69+200\cos^2t+20\cos t\sqrt{100\cos^2t-75}}.
\end{equation}
In order to provide the reader with a better opportunity to visualize the spatial 
form of the bipyramid $p(t)$ for small values of $t$, we expand the right-hand side of 
the formula (\ref{eqn:1}) in the Maclaurin series
\begin{equation}\label{eqn:2}
|AE|=3\sqrt{41}-\frac{117}{2\sqrt{41}}t^2+O(t^4),
\end{equation}
and note that $3\sqrt{41}\approx 19.209\ldots$.

We conclude this section with the computation of the volume of $p(t)$,
vol\,$p(t)$, as a function of $t$.
Obviously, area$(ABC)=\tfrac{1}{2}|AB||AC|\sin t$, where 
area$(ABC)$ stands for the area of the triangle $ABC$.
On the other hand, vol\,$p(t)=\tfrac{2}{3}|EF|\textrm{area}(ABC)$.
After simplifications, this yields
\begin{equation}\label{eqn:3}
\textrm{vol}\,p(t)=80(10\cos t+\sqrt{100\cos^2 t-75})\sin t.
\end{equation}
Expanding the right-hand side of the formula (\ref{eqn:3}) in the Maclaurin series,
we get
$$
\textrm{vol}\,p(t)=1200t-1400t^3+O(t^5).
$$
In particular, it follows from the latter formula that
vol\,$p(t)\to 0$ as $t\to 0$.

\textbf{4. Auxiliary polyhedron $q(t)$.}
In this section, we build and study a continuous family of auxiliary 
nonconvex oblique triangular bipyramids $q(t)$. 
The parameter $t$ is a real number, which takes values from
an interval $0<t<t'_0$.
The exact value of $t'_0$ will be discussed below.
\begin{figure}
\begin{center}
\includegraphics[width=0.4\textwidth]{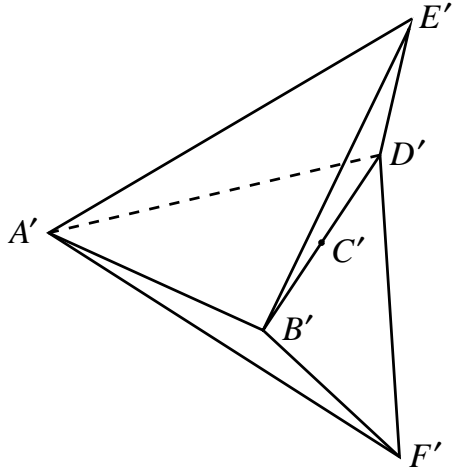}
\end{center}
\caption{Nonconvex oblique triangular bipyramid $q(t)$}
\end{figure}

The nonconvex oblique triangular bipyramid $q(t)$ is shown schematically in Fig.~2.
The point $C'$ is not a vertex of $q(t)$. $C'$ is the middle point of the edge $B'D'$.
It is shown in Fig.~2 for convenience of further explanations.
The bipyramid $q(t)$ is assumed to be symmetric with respect to the planes of the triangles
$A'B'D'$ and $A'E'F'$.
The lengths of the edges of $q(t)$ are defined in the following way:

(i) let $X$ and $Y$ be two points from the list $A,B,C,D,E,F$ 
(this is the list of the vertices of $p(t)$ extended by the point $C$);

(ii) suppose the straight line segment $XY$ entirely lies in $p(t)$
(i.\,e., $XY$ is either an edge of $p(t)$ or one of the following straight 
line segments: $BC$, $CD$, $CE$, $CF$);

(iii) suppose also that the straight line segment $X'Y'$ entirely lies in $q(t)$;

(iv) then we put by definition $|X'Y'|=|XY|$.

It follows immediately from this definition that $q(t)$ is intrinsically
isometric to $p(t)$ as soon as $p(t)$ and $q(t)$ exist.
From Sect.~3 we know that $p(t)$ does exist for $0\leq t<t_0=\pi/6$.
On the other hand, it is easy to observe that $q(t)$ exists for all $t$ such that 
$$
|C'E'|-|A'C'|<|A'E'|<|C'E'|+|A'C'|
$$ 
or, equivalently,
\begin{equation}\label{eqn:4}
3.3397\ldots\approx 12-5\sqrt{3}<|A'E'|<12+5\sqrt{3}\approx 20.6603\ldots.
\end{equation}
Using the formula (\ref{eqn:2}), we find $|A'E'|=3\sqrt{41}\approx 19.209\ldots$  for $t=0$ and conclude that (\ref{eqn:4}) holds true for  $t=0$.
It follows from (\ref{eqn:1}) that $|A'E'|$ is a continuous function in $t$ in a 
neighborhood of the point $t=0$. Thus, there exists a constant $t'_0>0$
such that the inequalities (\ref{eqn:4}) hold true for all $0<t<t'_0$. 
Hence, $q(t)$ does exist for $0\leq t<t'_0$.

We conclude this section with the computation of the volume of $q(t)$,
vol\,$q(t)$, as a function of $t$.
Obviously, 
\begin{equation}\label{eqn:5}
\text{\rm vol\,}q(t)=\frac{1}{3}|E'F'|\textrm{area}(A'B'D')=
\frac{2\sin\alpha}{3}|A'E'|\textrm{area}(A'B'D'),
\end{equation}
where $\alpha$ is the angle of the triangle $A'C'E'$ at the vertex $A'$.
Since $A'B'D'$ is a regular triangle with the side length equal to 10,
we have $\textrm{area}(A'B'D')=25\sqrt{3}$ 
is independent of $t$).
Using the cosine law for the triangle $A'C'E'$, we find 
\begin{equation*}
\cos\alpha=\frac{|A'E'|^2-69}{10\sqrt{3}|A'E'|}.\quad \textrm{Thus}\quad 
\sin\alpha=\frac{\sqrt{438|A'E'|^2-|A'E'|^4-69^2}}{10\sqrt{3}|A'E'|}.
\end{equation*}
Substituting the expressions for $\textrm{area}(A'B'D')$ and $\sin\alpha$ 
to (\ref{eqn:5}),
we get 
\begin{equation}\label{eqn:6}
\textrm{vol}\,q(t)=\frac{5}{3}\sqrt{438|A'E'|^2-|A'E'|^4-69^2}.
\end{equation}
Recall that $|A'E'|=|AE|$ and $|AE|$ is the function in $t$ defined by the formula
(\ref{eqn:1}). 
Substitute (\ref{eqn:1}) in (\ref{eqn:6}) and expand the right-hand side of the 
formula obtained in the Maclaurin series (the result of substitution itself is too 
complicated and we prefer not to wright it here).
The first two members of the Maclaurin series are as follows
$$
\textrm{vol}\,q(t)=50\sqrt{23}+\frac{3750}{\sqrt{23}}t^2+O(t^4).
$$
In particular, we conclude from the latter formula that 
vol\,$q(t)\to 50\sqrt{23}\approx 239.79$ as $t\to 0$.

\textbf{6. Proof of Theorem~1.}
Given a constant $c>0$, find $t^*$ such that
$0<t^*<\min\{t_0,t_0'\}$ and 
\begin{equation}\label{eqn:7}
\frac{\text{\rm vol\,}q(t^*)}{\text{\rm vol\,}p(t^*)}>c.
\end{equation}
We can satisfy (\ref{eqn:7}) because we know from Sect.~4 that 
vol\,$q(t)\to 50\sqrt{23}$ as $t\to 0$
and we know from Sect.~3 that vol\,$p(t)\to 0$ as $t\to 0$.
From Sect.~4 we also know that $p(t^*)$ and $q(t^*)$ are intrinsically isometric to
each other. By definition, put $P=p(t^*)$ and $Q=q(t^*)$.
The bipyramids $P$ and $Q$ satisfy all the conditions of Theorem~1. \hfill $\square$

\bigskip

\noindent{Victor Alexandrov}

\noindent\textit{Sobolev Institute of Mathematics}

\noindent\textit{Koptyug ave., 4}

\noindent\textit{Novosibirsk, 630090, Russia}

and

\noindent\textit{Department of Physics}

\noindent\textit{Novosibirsk State University}

\noindent\textit{Pirogov st., 2}

\noindent\textit{Novosibirsk, 630090, Russia}

\noindent\textit{e-mail: alex@math.nsc.ru}

\bigskip

\noindent{Originally submitted to arXive: July 22, 2016}

\end{document}